\documentclass[12pt]{article}
\usepackage{pdfpages}
\usepackage{url}
\usepackage{cite}
\usepackage[margin=2.5cm]{geometry}
\usepackage[normalem]{ulem}
\usepackage[T1]{fontenc}
\usepackage[utf8]{inputenc}
\usepackage[english]{babel}
\usepackage{url}
\usepackage{graphicx}
\usepackage{mathtools}
\usepackage{amsmath}
\usepackage{amssymb}
\usepackage{amsthm}
\usepackage{dsfont}
\usepackage{amsfonts}
\usepackage{mathrsfs}
\usepackage{bbm}
\usepackage{enumerate}
\usepackage{graphicx,tikz,pgf}
\usepackage[colorlinks = true,
            linkcolor = red,
            urlcolor  = red,
            citecolor = red,
            anchorcolor = red]{hyperref}

\usetikzlibrary{arrows}

\usepackage[all,cmtip]{xy}
\usepackage{amssymb,amsmath,amsthm,mathrsfs}
\def\RR{{\mathbb R}}

\def\1{{\mathbbm 1}}

\def\diff{{\rm Diff}_+}

\def\diffs1{\diff(S^1)}

\def\psl2r{{\rm PSL}(2,\RR)}
\def\sl2r{{\rm SL}(2,\RR)}
\def\su11{{\rm SU}(1,1)}
\def\2dmob{{\overline{\psl2r}\times\overline{\psl2r}}}
\def\<{\langle}
\def\>{\rangle}

\theoremstyle{definition}

\theoremstyle{remark}

\theoremstyle{definition}

\begin{document}

\title{A correction of the historiographical record on the probability integral
}

\markright{The probability integral}
\author{
{\bf Fausto Di Biase}\footnote{Supported in part by  INdAM-GNAMPA and by Fondi di Ricerca di Ateneo Universit\`a ``G.\! D'Annunzio''},\\
Dipartimento di Economia,
Universit\`a ``G.\! D'Annunzio'' di Chieti-Pescara,\\
Viale Pindaro, 42, I-65127, Pescara, Italy\\
email: {\tt fdibiase@unich.it}\\
}

\maketitle

\begin{abstract}
We correct a common (but mistaken) attribution of the evaluation of 
the probability integral, usually attributed to Poisson, Gauss, or 
Laplace.
\footnote{MSC: 01A85; 01A50 } 
\end{abstract}

\section{Common but mistaken attributions}

The  \textit{probability integral} is the following one:
\begin{equation}
\int_\mathbb{R} e^{-x^2}\,dx=\sqrt{\pi}
\tag{E}
\end{equation}
Some modern authors attribute the probability integral~(E) to \textit{Poisson} \cite[p. 132]{Gnedenko}. 
However, Gauss himself (in an important book published in 1809, to which we shall return momentarily) attributes the result to 
Laplace. 
The work of Laplace was published in 1774 \cite[pp.35-36]{Laplace_1774}.  In 1812 Laplace returned to this calculation and gave a different proof \cite[pp.93--96]{Laplace_1812}. Even modern authors attribute the result to 
Laplace. For example, S. Stigler, in his study of Laplace's work,   
writes 

\begin{quotation}
``Since the proof includes what seems to be the earliest evaluation of the definite integral 
$
2\int_0^{\infty}\frac{1}{\sqrt{2\pi}\sigma}\exp\left(
-\frac{z^2}{2\sigma^2}
\right)
dz
=
1
$ 
Laplace may have been the first to integrate the normal density'' \cite[p.360]{Stigler_1986}.
\end{quotation}
Currently a consensus points to Laplace as the author of the evaluation of~(E). We will show that attribution of~(E) to Laplace given by Gauss is mistaken. 

\section{The correction}

Firstly, observe that Laplace~\cite[p.36]{Laplace_1774} proves~(E) in the form~(E.0)
\begin{equation}
\int_0^1
{
[\ln
(
{1/x}
]
}
^
{-1/2}
dx
=
\sqrt{\pi}
\tag{E.0}
\end{equation}
Indeed, the substitution $x=e^{-t^2}$ shows that~(E.0) is equivalent to~(E). 

At the beginning of his proof, Laplace gives a generic attribution to Euler, and writes 
``{Voir le \textit{Calcul int{\'e}gral} de M. Euler}'' \cite[p.35]{Laplace_1774}. 

As a matter of fact, in \cite[p. 111]{Euler_1772}, published in 1772, \textit{Euler had already proved}~(E.0). Moreover, in his previous work~\cite{Euler_1738}, where he introduced the factorial of non integer numbers, now known as the Gamma function, 
Euler had 
proved that 
\begin{equation}
\int_0^1
{
[
\ln(1/x)
]
}^{1/2}
dx
=\frac{1}{2}\sqrt{\pi}
\tag{E.0'}
\end{equation}
which is again equivalent to~(E).\footnote{Indeed, the Gamma function may be written as 
$
\Gamma(z)=\int_0^1
{
[
\ln(1/t)
]
}^{z-1}
dt
$
and therefore the expression in~(E.0') is equal to $\Gamma(3/2)$ which is 
equal to $\frac{1}{2}\Gamma(1/2)=\sqrt{\pi}/2$. All of this had been proved by Euler, although in a different notation.}

%
%

\subsection{Who knew it and when?}

Mathematicians who were contemporaries of Gauss were very well aware of the fact that~(E.0) is equivalent to~(E). Indeed, inside a printed copy of Gauss's book \cite{Gauss_1809}, located in the library of the 
ETH Z{\"u}rich, a handwritten footnote by Barnaba Oriani, a contemporary of Gauss, written in Latin, states  that the evaluation in~(E) is really due to Euler, not to Laplace. 
This handwritten footnote is located on page 212, where Gauss 
mentions this result as \textit{theorema elegans primo ab ill. Laplace inventum} (The elegant theorem first proved by illustrious Laplace). 
Here follows a translation of Oriani's footnote. 

\begin{quotation}
\textit{
The elegant theorem, attributed to the illustrious La Place, was in reality first found by Leonard Euler. As a matter of fact, Euler had proved it first, in Comment. Acad. Petrop. (vol 16) that 
$\displaystyle{
{-\int \frac{1}{\sqrt{\ln{\frac{1}{x}}}}}
\,dx
}$
extended from $x=1$ to $x=0$ is equal to $\pi$, i.e., half of the length of a circumference of radius one. Hence, putting  $x=e^{-t^2}$ one has 
$\frac{
-dx
}
{
\sqrt{
\ln{\frac{1}{x}}
}
}
=2e^{-t^2}dt
$. 
Hence the integral $\int e^{-t^2}$ from $t=0$ to $t=\infty$ is 
$=\frac{1}{2}\sqrt{\pi}$ and therefore the same integral from $t=-\infty$ to $t=+\infty$ will be $=\sqrt{\pi}$.
}
\end{quotation}

Barnaba Oriani (1752-1832) was a mathematician and astronomer. He was a friend of Giuseppe Piazzi (1746,1826), the astronomer who on January 1, 1801, 
from his observatory in Palermo,  
discovered the planetoid Ceres, located between the orbits of Mars and Jupiter, and collected the data of its position until February 11, when ``because of sickness, it was not possible for Piazzi to continue observations'' \cite[p.25]{Herrmann}. ``When the news of the discovery reached northern Europe a period of bad weather hindered further monitoring of the object; and finally when there was another opportunity to make observations the new object could not be found'' \cite[p.25]{Herrmann}. On the basis of few data, Gauss was able to forecast the position of Ceres. ``This led to its rediscovery 
by Olbers in Bremen on January 1, 1802, precisely a year after the first discovery [...]'' \cite[p.26]{Herrmann}. ``The astronomers were exceedingly impressed by this achievement, for the ellipse calculated by Gauss represented the observations of Piazzi astonishingly well. For the time being Gauss did not communicate details of his new method, but everyone was quite convinced that a new method had been developed.'' \cite[p.26]{Herrmann}. ``The secret of his method was soon to be disclosed: in 1809 Gauss published his procedure in his classical work \textit{Theoria Motus Corporum Coelestium in sectionibus conicis solem ambientium} (Theory of the Motion of Celestial bodies which revolve  about the Sun in Conic Sections).\cite[p.26]{Herrmann}. 
This is the book where Gauss attributes the calculation in~(E) to 
Laplace.

\bibliography{GaussianBib}{}
\bibliographystyle{apalike}

\end{document}